\documentstyle{article}
\begin{document}
 \large
\begin{center}{ON CLASSIFICATION OF FINITE DIMENSIONAL COMPLEX FILIFORM LEIBNIZ
ALGEBRAS (PartI) } \end{center}
\begin{center} {Bekbaev U.D.\footnote[1]{e-mail: bekbaev@science.upm.edu.my},
Rakhimov I.S.\footnote[2]{e-mail: isamiddin@science.upm.edu.my}}
\end{center}
\begin{center} {Department of Mathematics $\&$ Institute for Mathematical
Research,}\end{center}
\begin{center} {FS,UPM, 43400, Serdang, Selangor Darul Ehsan, Malaysia.}\end{center}

\begin{abstract}
The paper is devoted to classification problem of finite
dimensional complex none Lie filiform Leibniz algebras. The
motivation to write this paper is an unpublished yet result of
J.R.Gomez, B.A.Omirov [1] on necessary and sufficient conditions
for two finite dimensional complex filiform Leibniz algebras to be
isomorphic. We suggest another approach to this problem. The
approach that we use for classification is in terms of invariants.
In fact, utilizing this method for any given low dimensional case
all filiform Leibniz algebras can be classified. Moreover, the
results can be used for geometrical classification of orbits of
such algebras.
 \end{abstract}



\section{Introduction}
This paper aims to investigate a class of nonassociative algebras
which generalizes the class of Lie algebras. These algebras
satisfy certain identities that were suggested by J.-L.Loday [2].
When he used the tensor product instead of external product in the
definition of the $n$-th cochain, in order to prove the
differential property, that is defined on cochains, it sufficed to
replace the anticommutativity and Jacoby identity by the Leibniz
identity. This is an essential one of the motivation to appear for
this class of algebras.

In this paper we suggest an algebraic approach to the
classification problem for filiform Leibniz algebras. Utilizing
this method for any fixed low dimensional case the corresponding
classes of filiform Leibniz algebras can be classified completely.
Moreover, the results may be used for geometric classification in
the sense of geometric invariant theory [3]. It is assumed that it
will be the subject of one of the next papers. For geometric
classification of complex nilpotent Leibniz algebras of dimension
at most four we refer to [4].

Let $V$ be a vector space of dimension $n$ over an algebraically
closed field $K$ (char$K$=0). The bilinear maps $V \times V
\rightarrow V$ form a vector space $Hom(V\otimes V,V)$ of
dimension $n^{3}$, which can be considered together with its
natural structure of an affine algebraic variety over $K$ and
denoted by $Alg_{n}(K)\cong K^{n^{3}}$. An $n$-dimensional algebra
$L$ over $K$ may be considered as an element $\lambda(L)$ of
$Alg_{n}(K)$ via the bilinear mapping $\lambda:L\otimes L\to L$
defining an binary algebraic operation on $L:$ let $\{
e_1,e_2,\ldots,e_n\}$ be a basis of the algebra $L.$ Then the
table of multiplication of $L$ is represented by point
$(\gamma_{ij}^{k})$ of this affine space as follow:
$$\lambda(e_i,e_j)=\sum\limits_{k=1}^n\gamma_{ij}^ke_k.$$
$\gamma_{ij}^k$ are called {\it structural constants} of $L.$ The
linear reductive group $GL_n(K)$ acts on $Alg_{n}(K)$ by
$(g*\lambda)(x,y)=g(\lambda(g^{-1}(x),g^{-1}(y)))$("transport of
structure"). Two algebras $\lambda_{1}$ and $\lambda_{2}$ are
isomorphic if and only if they belong to the same orbit under this
action. The orbit of $\lambda$ under this action is denoted by
$O(\lambda)$. It is clear that elements of the given orbit are
isomorphic to each other algebras. The classification means to
specify the representatives of the orbits. A simple criterion, to
decide if the given two algebras are isomorphic, is desired.

\section{Preliminaries}

{\bf Definition 1.} An algebra $L$ over a field $K$ is called a
{\it Leibniz algebra} if it satisfies the following Leibniz
identity: $$[x,[y,z]]=[[x,y],z]-[[x,z],y],$$ where $[\cdot,\cdot]$
denotes the multiplication in $L.$ Let $Leib_{n}(K)$ be a
subvariety of $Alg_{n}(K)$ consisting of all $n$-dimensional
Leibniz algebras over K. It is invariant under the above mentioned
action of $GL_n(K)$. As a subset of $Alg_{n}(K)$ the set
$Leib_{n}(K)$ is specified by system of equations with respect to
structural constants $\gamma_{ij}^{k}$:
$$\sum\limits_{\emph{l}=1}^{\emph{n}}{(\gamma_{\emph{jk}}^{\emph{l}}
\gamma_{\emph{il}}^{\emph{m}}-\gamma_{\emph{ij}}^{\emph{l}}\gamma_{\emph{lk}}^{\emph{m}}+
\gamma_{\emph{ik}}^{\emph{l}}\gamma_{\emph{lj}}^{\emph{m}})}=0$$
It is easy to see that if the bracket in Leibniz algebra happens
to be anticommutative then it is a Lie algebra. So Leibniz
algebras are "noncommutative" generalization of Lie algebras. As
to classifications of low dimensional Lie algebras they are well
known. But unless simple Lie algebras the classification problem
of all Lie algebras in common remains a big problem. Yu.I.Malcev
[4] reduced the classification of solvable Lie algebras to the
classification of nilpotent Lie algebras. Apparently the first
non-trivial classification of some classes of low-dimensional
nilpotent Lie algebra are due to Umlauf. In his thesis [6] he
presented the redundant list of nilpotent Lie algebras of
dimension at most seven. He gave also the list of nilpotent Lie
algebras of dimension less than ten admitting so-called adapted
basis (now, the nilpotent Lie algebras with this property are
called \emph{filiform Lie algebras}). It was shown by M.Vergne [7]
the impotentness of filiform Lie algebras in the study of variety
of nilpotent Lie algebras laws. Up to now the several
classifications of low-dimensional nilpotent Lie algebras have
been done. Unfortunately, many of these papers are based on direct
computations (by hand) and the complexity of those computations
leads frequently to errors. We refer the reader to [8] for
comments and corrections of the classification errors.

Further if it is not asserted additionally all algebras assumed to
be over the field of complex numbers.

Let $L$ be a Leibniz algebra. We put: $$ L^1=L,\quad
L^{k+1}=[L^k,L],\enskip k\in N.$$

{\bf Definition 2.} A Leibniz algebra $L$ is said to be {\it
nilpotent} if there exists an integer $s\in N,$ such that $L^{1}
\supset L^{2} \supset ... \supset L^{s}=\{0\}.$ The smallest
integer $s$ for which $L^{s}=0$ is called {\it the nilindex} of
$L$.

{\bf Definition 3.} An $n$-dimensional Leibniz algebra $L$ is said
to be {\it filiform} if $dim L^i =n-i,$ where $2\le i\le n.$

{\bf Theorem 1[9].} Arbitrary complex non-Lie filiform Leibniz
algebra of dimension $n+1$ is isomorphic to one of the following
filiform Leibniz algebras:
$$
\mbox{a) (The first class):} \ \ \left\{
\begin{array}{lll}
[e_{0},e_{0}]=e_{2},\\

[e_i,e_{0}]=e_{i+1}, \qquad \qquad \qquad \qquad \qquad
\qquad \qquad  \ \ 1\leq i\leq n-1 \\

[e_{0},e_{1}]= \alpha_{3} e_{3}+ \alpha_{4} e_{4}+...+\alpha_{n-1}
e_{n-1}+
\theta e_{n}, \\

[e_{j},e_{1}]=\alpha_{3}e_{j+2}+
\alpha_{4}e_{j+3}+...+\alpha_{n+1-j}e_{n}, \qquad 1\leq j\leq n-2

\end{array}
\right.
$$
$$
\mbox{b) (The second class):} \ \ \left\{
\begin{array}{lll}
[\emph{e}_{0},\emph{e}_{0}]=\emph{e}_{2},\\

[e_{i},e_{0}]=e_{i+1}, \qquad \qquad \qquad \qquad \qquad
\qquad \qquad  \ \ 2\leq i\leq n-1 \\

[e_{0},e_{1}]= \beta_{3} e_{3}+ \beta_{4}
e_{4}+...+\beta_{n} e_{n}, \\

[e_{1},e_{1}]=\gamma e_{n},\\

[e_{j},e_{1}]=\beta_{3}e_{j+2}+
\beta_{4}e_{j+3}+...+\beta_{n+1-j}e_{n}, \qquad 2\leq j\leq n-2

\end{array}
\right.
$$
where $[,]$ is the multiplication in a Leibniz algebra, $\{
e_0,e_1, \ldots, e_n\}$ is a basis of the algebra and omitted
products are assumed to be zero.

Note that the algebras from the first class and the second class
never are isomorphic to each other.

In this paper we will deal with the first class of algebras, as
for the second class it will be considered in our another paper.

Let us denote by $L(\alpha) $,  the $(n+1)$-dimensional filiform
non-Lie Leibniz algebra defined by parameters
$\alpha=(\alpha_{3},\alpha_{4},...,\alpha_{\emph{n}},\theta)$. The
set of all $(n+1)$-dimensional complex filiform Leibniz algebras
is denoted by $FL_{n+1}$. It is a closed and invariant subset of
the variety of nilpotent Leibniz algebras.

Using the method of simplification of the basis transformations in
[1] the following criterion on isomorphism of two
$(n+1)$-dimensional filiform Leibniz algebras was given. Namely:
let $n\geq 3$.

{\bf Theorem 2.}[1] Two algebras $L(\alpha)$ and $L(\alpha')$ from
$FL_{n+1}$, where
$\alpha=(\alpha_{3},\alpha_{4},...,\alpha_{\emph{n}},\theta)$ and
$\alpha'=(\alpha'_{3},\alpha'_{4},...,\alpha'_{\emph{n}},\theta')$,
are isomorphic if and only if there exist complex numbers
$\emph{A},\emph{B}$ such that $\emph{A}(\emph{A}+\emph{B})\neq$ 0
and the following conditions hold:
$$
\left\{
\begin{array}{lll}

\alpha'_{3}=\frac{A+B}{A^{2}}\alpha_{3},\\
\alpha'_{t}=\frac{1}{A^{t-1}}((A+B)\alpha_{t}- \sum
\limits_{k=3}^{t-1}(C_{k-1}^{k-2}A^{k-2}B \alpha_{t+2-k}+
C_{k-1}^{k-3}A^{k-3}B^{2} \sum \limits_{i_{1}=k+2}^{t}
\alpha_{t+3-i_{1}}\cdot\alpha_{i_{1}+1-k}+\\
C_{k-1}^{k-4}A^{k-4}B^{3}\sum\limits_{i_{2}=k+3}^{t}
\sum\limits_{i_{1}=k+3}^{i_{2}}\alpha_{t+3-i_{2}}\cdot\alpha_{i_{2}+3-i_{1}}\cdot
\alpha_{i_{1}-k}+...+\\
C_{k-1}^{1}AB^{k-2}\sum\limits_{i_{k-3}=2k-2}^{t}
\sum\limits_{i_{k-4}=2k-2}^{i_{k-3}}...
\sum\limits_{i_{1}=2k-2}^{i_{2}}\alpha_{t+3-i_{k-3}}\cdot\alpha_{i_{k-3}+3-i_{k-4}}
\cdot...\cdot\alpha_{i_{2}+3-i_{1}}\cdot\alpha_{i_{1}+5-2k}\\+B^{k-1}
\sum\limits_{i_{k-2}=2k-1}^{t}
\sum\limits_{i_{k-3}=2k-1}^{i_{k-2}}...
\sum\limits_{i_{1}=2k-1}^{i_{2}}\alpha_{t+3-i_{k-2}}\cdot\alpha_{i_{k-2}+3-i_{k-3}}
\cdot...\cdot\alpha_{i_{2}+3-i_{1}}\alpha_{i_{1}+4-2k})\cdot\alpha'_{k}), \\

\mbox{where} \ \ 4 \leq t \leq n.\\

\theta'=\frac{1}{A^{n-1}}(A\theta+B\alpha_{n}- \sum
\limits_{k=3}^{n-1}(C_{k-1}^{k-2}A^{k-2}B \alpha_{n+2-k}+
C_{k-1}^{k-3}A^{k-3}B^{2} \sum \limits_{i_{1}=k+2}^{n}
\alpha_{n+3-i_{1}}\cdot\alpha_{i_{1}+1-k}+\\
C_{k-1}^{k-4}A^{k-4}B^{3}\sum\limits_{i_{2}=k+3}^{n}
\sum\limits_{i_{1}=k+3}^{i_{2}}\alpha_{n+3-i_{2}}\cdot\alpha_{i_{2}+3-i_{1}}
\cdot\alpha_{i_{1}-k}+...+\\
C_{k-1}^{1}AB^{k-2}\sum\limits_{i_{k-3}=2k-2}^{n}
\sum\limits_{i_{k-4}=2k-2}^{i_{k-3}}...
\sum\limits_{i_{1}=2k-2}^{i_{2}}\alpha_{n+3-i_{k-3}}\cdot\alpha_{i_{k-3}+3-i_{k-4}}
\cdot...\cdot\alpha_{i_{2}+3-i_{1}}\cdot\alpha_{i_{1}+5-2k}\\+B^{k-1}
\sum\limits_{i_{k-2}=2k-1}^{n}
\sum\limits_{i_{k-3}=2k-1}^{i_{k-2}}...
\sum\limits_{i_{1}=2k-1}^{i_{2}}\alpha_{n+3-i_{k-2}}\cdot\alpha_{i_{k-2}+3-i_{k-3}}
\cdot...\cdot\alpha_{i_{2}+3-i_{1}}\cdot\alpha_{i_{1}+4-2k})\cdot\alpha'_{k}),

\end{array}
\right.
$$

Here are the above systems of equalities for some low dimensional
cases:

Case of $n=4$ i.e. $\textbf{dimL=5:}$
$$ \begin{array}{lll} \left\{
\begin{array}{lll}
\alpha'_{3}=\frac{1}{A}(1+\frac{B}{A})\alpha_{3},\\
\\
\alpha'_{4}=\frac{1}{A^{2}}(1+\frac{B}{A})(\alpha_{4}-2\frac{B}{A}\alpha_{3}^{2})\\
\\
\theta'=\frac{1}{A^2}[\theta+\frac{B}{A}\alpha_{4}-2(1+\frac{B}{A})\frac{B}{A}\alpha_{3}^{2}].
\end{array}
\right.
\end{array}
$$
Case of $n=5$ i.e. $\textbf{dimL=6:}$
$$ \begin{array}{lll} \left\{
\begin{array}{lll}
\alpha'_{3}=\frac{1}{A}(1+\frac{B}{A})\alpha_{3},\\
\\
\alpha'_{4}=\frac{1}{A^{2}}(1+\frac{B}{A})(\alpha_{4}-2\frac{B}{A}\alpha_{3}^{2})\\
\\
\alpha'_{5}=\frac{1}{A^{3}}(1+\frac{B}{A})[\alpha_{5}-5\frac{B}{A}
(\alpha_{4}-\frac{B}{A}\alpha_{3}^{2})\alpha_{3}]\\
\\
\theta'=\frac{1}{A^3}[\theta+\frac{B}{A}\alpha_{5}-5(1+\frac{B}{A})
\frac{B}{A}(\alpha_{4}-\frac{B}{A}\alpha_{3}^{2})\alpha_{3}].
\end{array}
\right.
\end{array}
$$

Case of $n=6$ i.e. $\textbf{dimL=7:}$
$$ \begin{array}{lll} \left\{
\begin{array}{lll}
\alpha'_{3}=\frac{1}{A}(1+\frac{B}{A})\alpha_{3},\\
\\
\alpha'_{4}=\frac{1}{A^{2}}(1+\frac{B}{A})(\alpha_{4}-2\frac{B}{A}\alpha_{3}^{2})\\
\\
\alpha'_{5}=\frac{1}{A^{3}}(1+\frac{B}{A})[\alpha_{5}-5\frac{B}{A}
(\alpha_{4}-\frac{B}{A}\alpha_{3}^{2})\alpha_{3}]\\
\\
\alpha'_{6}=\frac{1}{A^{4}}(1+\frac{B}{A})[\alpha_{6}-6\frac{B}{A}\alpha_3\alpha_5+
21(\frac{B}{A})^2\alpha_3^2\alpha_4-3\frac{B}{A}\alpha_4^2-14(\frac{B}{A})^3\alpha_{3}^{4}]\\
\\
\theta'=\frac{1}{A^4}\{\theta+\frac{B}{A}\alpha_{6}-(1+\frac{B}{A})
[6\frac{B}{A}\alpha_3\alpha_5-21(\frac{B}{A})^2\alpha_3^2\alpha_4+3\frac{B}{A}\alpha_4^2
+14(\frac{B}{A})^3\alpha_{3}^{4}]\}.
\end{array}
\right.
\end{array}
$$

Case of $n=7$ i.e. $\textbf{dimL=8:}$
$$ \begin{array}{llr} \left\{
\begin{array}{lll}
\alpha'_{3}=\frac{1}{A}(1+\frac{B}{A})\alpha_{3},\\
\\
\alpha'_{4}=\frac{1}{A^{2}}(1+\frac{B}{A})(\alpha_{4}-2\frac{B}{A}\alpha_{3}^{2})\\
\\
\alpha'_{5}=\frac{1}{A^{3}}(1+\frac{B}{A})[\alpha_{5}-5\frac{B}{A}
(\alpha_{4}-\frac{B}{A}\alpha_{3}^{2})\alpha_{3}]\\
\\
\alpha'_{6}=\frac{1}{A^{4}}(1+\frac{B}{A})[\alpha_{6}-6\frac{B}{A}\alpha_3\alpha_5+
21(\frac{B}{A})^2\alpha_3^2\alpha_4-3\frac{B}{A}\alpha_4^2-14(\frac{B}{A})^3\alpha_{3}^{4}]\\
\\
\alpha'_{7}=\frac{1}{A^{5}}(1+\frac{B}{A})[\alpha_{7}-7\frac{B}{A}\alpha_3\alpha_6+
28(\frac{B}{A})^2\alpha_3^2\alpha_5+28(\frac{B}{A})^2\alpha_3\alpha_4^2-
7\frac{B}{A}\alpha_4\alpha_5\\
-84(\frac{B}{A})^3\alpha_{3}^{3}\alpha_4+42(\frac{B}{A})^4\alpha_{3}^{5}]\\
\\
\theta'=\frac{1}{A^5}
\{\theta+\frac{B}{A}\alpha_{7}-(1+\frac{B}{A})
[7\frac{B}{A}\alpha_3\alpha_6-28(\frac{B}{A})^2\alpha_3^2\alpha_5-28(\frac{B}{A})^2\alpha_3\alpha_4^2+
7\frac{B}{A}\alpha_4\alpha_5\\
+84(\frac{B}{A})^3\alpha_{3}^{3}\alpha_4-42(\frac{B}{A})^4\alpha_{3}^{5}]\}.
\end{array}
\right.
\end{array}
$$

To deal with the classification of $FL_{n+1}$ with respect to the
above mentioned action we represent it as a disjoint union of an
open and closed (with respect to the Zarisski topology) subsets.
Moreover each of these subsets are invariant under the
corresponding transformations presented in Theorem 2. Then we
formulate the solution of the isomorphism problem for the
corresponding algebras from the open subset. Similar approach can
be used to solve isomorphism problem for the algebras from the
corresponding closed subset.

It is not difficult to notice that the expressions for
$\alpha'_{t}$, $\theta'$ in Theorem 2 can be represented in the
following form:
$$\alpha'_{t}=\frac{1}{A^{t-2}}\varphi_{t}(\frac{B}{A};\alpha),$$ where
$\alpha=(\alpha_{3},\alpha_{4},...,\alpha_{n}, \theta)$ and
$\varphi_{t}(y;z)=\varphi_{t}(y;z_{3},z_{4},...,z_{n},z_{n+1})=$
$$\begin{array}{lll}
((1+y)z_{t}- \sum \limits_{k=3}^{t-1}(C_{k-1}^{k-2}y z_{t+2-k}+
C_{k-1}^{k-3}y^{2} \sum \limits_{i_{1}=k+2}^{t} z_{t+3-i_{1}}
\cdot z_{i_{1}+1-k}+\\
\\
C_{k-1}^{k-4}y^{3}\sum\limits_{i_{2}=k+3}^{t}
\sum\limits_{i_{1}=k+3}^{i_{2}} z_{t+3-i_{2}} \cdot
z_{i_{2}+3-i_{1}} \cdot z_{i_{1}-k}+...+\\
\\
C_{k-1}^{1}y^{k-2}\sum\limits_{i_{k-3}=2k-2}^{t}
\sum\limits_{i_{k-4}=2k-2}^{i_{k-3}}...
\sum\limits_{i_{1}=2k-2}^{i_{2}} z_{t+3-i_{k-3}} \cdot
z_{i_{k-3}+3-i_{k-4}} \cdot...\cdot z_{i_{2}+3-i_{1}} \cdot
z_{i_{1}+5-2k}+\\
\\y^{k-1} \sum\limits_{i_{k-2}=2k-1}^{t}
\sum\limits_{i_{k-3}=2k-1}^{i_{k-2}}...
\sum\limits_{i_{1}=2k-1}^{i_{2}} z_{t+3-i_{k-2}} \cdot
z_{i_{k-2}+3-i_{k-3}} \cdot...\cdot z_{i_{2}+3-i_{1}}
z_{i_{1}+4-2k}) \cdot \varphi_{k}(y;z)),\\
\\
\mbox{for} \ 3 \leq t \leq n.\end{array}$$
$\theta'=\frac{1}{A^{n-2}}\varphi_{n+1}(\frac{B}{A};\alpha),
\mbox{where}\
\varphi_{n+1}(y;z)=\varphi_{n+1}(y;z_{3},z_{4},...,z_{n},z_{n+1})=$
$$\begin{array}{lll}
(z_{n+1}+y z_{n}- (1+y)\sum \limits_{k=3}^{n-1}(C_{k-1}^{k-2}y
z_{n+2-k}+ C_{k-1}^{k-3}y^{2} \sum \limits_{i_{1}=k+2}^{n}
z_{n+3-i_{1}} \cdot z_{i_{1}+1-k}+\\
\\
C_{k-1}^{k-4}y^{3}\sum\limits_{i_{2}=k+3}^{n}
\sum\limits_{i_{1}=k+3}^{i_{2}} z_{n+3-i_{2}} \cdot
z_{i_{2}+3-i_{1}}
\cdot z_{i_{1}-k}+...+\\
\\
C_{k-1}^{1}y^{k-2}\sum\limits_{i_{k-3}=2k-2}^{n}
\sum\limits_{i_{k-4}=2k-2}^{i_{k-3}}...
\sum\limits_{i_{1}=2k-2}^{i_{2}}z_{n+3-i_{k-3}} \cdot
z_{i_{k-3}+3-i_{k-4}} \cdot...\cdot z_{i_{2}+3-i_{1}}\cdot
z_{i_{1}+5-2k}\\
\\
+y^{k-1} \sum\limits_{i_{k-2}=2k-1}^{n}
\sum\limits_{i_{k-3}=2k-1}^{i_{k-2}}...
\sum\limits_{i_{1}=2k-1}^{i_{2}}z_{n+3-i_{k-2}}\cdot
z_{i_{k-2}+3-i_{k-3}} \cdot...\cdot z_{i_{2}+3-i_{1}}\cdot
z_{i_{1}+4-2k})\cdot\varphi_{k}(y;z)).\end{array}$$

To simplify notation let us agree that in the above case for
transition from the $(n+1)$-dimensional filiform Leibniz algebra
$L(\alpha)$ to the $(n+1)$-dimensional filiform Leibniz algebra
$L(\alpha')$ we write
$\alpha'=\rho(\frac{1}{A},\frac{B}{A};\alpha)$, \ \ where
$\alpha=(\alpha_{3},\alpha_{4},...,\alpha_{n},\theta)$,

$$\rho(\frac{1}{A},\frac{B}{A};\alpha)=(\rho_{1}(\frac{1}{A},\frac{B}{A};\alpha),\rho_{2}(\frac{1}{A},\frac{B}{A};\alpha),...,\rho_{n-1}(\frac{1}{A},\frac{B}{A};
\alpha)),$$

$$\rho_{t}(x,y;z)=x^{t}\varphi_{t+2}(y;z)\ \mbox{for} \ 1 \leq t \leq
n-2$$ \ and \

$$\rho_{n-1}(x,y;z)=x^{n-2}\varphi_{n+1}(y,z)$$

Here are the main properties of the operator $\rho$, derived from
the fact that $\rho(\frac{1}{A},\frac{B}{A};\cdot)$ is an action
of a group, that will be used later on.

$$\begin{array}{lll}{1^{0}. \ \ \rho(1,0;\cdot) \ \  \mbox{is the identity
operator}.}\\
\\ 2^{0}. \
\rho(\frac{1}{A_{2}},\frac{B_{2}}{A_{2}};\rho(\frac{1}{A_{1}},
\frac{B_{1}}{A_{1}};\alpha))=\rho(\frac{1}{A_{1}A_{2}},
\frac{A_{1}B_{2}+A_{2}B_{1}+B_{1}B_{2}}{A_{1}A_{2}};\alpha).\\
\\ 3^{0}. \ \ \mbox{If} \ \
\alpha'=\rho(\frac{1}{A},\frac{B}{A};\alpha) \ \ \mbox{then} \ \
\alpha=\rho(A,-\frac{B}{A+B};\alpha').\end{array}$$

\section{Classification theorems}

{\bf Definition 4.} An action of algebraic group $G$ on a variety
$Z$ is a morphism

$\sigma: G \times Z \longrightarrow Z$ with

$(i)\ \  \sigma(e,z)=z$, where $e$ is the unit element of $G$ and
$z \in Z$

$(ii)\ \ \sigma(g,\sigma(h,z))=\sigma(gh,z)$, for any $g,h \in G$
and $z \in Z$.

We shortly write $gz$ for $\sigma(g,z)$, and call $Z$ a
$G$-variety.\\

{\bf Definition 5.} A morphism $\emph{f}:Z \longrightarrow K$ is
said to be invariant if $\emph{f}(gz)= \emph{f}(z)$ for any $g \in
G$ and $z \in Z$.

The algebra of invariant morphisms on $Z$ with respect to the
action of the group G is denoted by $K[Z]^{G}$. Sometimes this
algebra is a finitely generated $K$-algebra. This is referred to
in [10] as the "first fundamental problem of invariant theory". If
$Z$ is an irreducible then the field of rational invariants can be
defined as a quotient field of $K[Z]^{G}$. It is always finitely
generated as a subalgebra of the finitely generated algebra
$K(Z)$. Description the field of rational invariants is an another
important classical problem of the invariant theory [11].

Actually, we use some elements of the algebra of invariant
morphisms under the above mentioned adapted action on the variety
of filiform Leibniz algebras to solve isomorphism problem.

From here on $n$ is a positive integer. We assume that $n\geq4$
since there are complete classifications of complex nilpotent
Leibniz algebras of dimension at most four [12],[13].

We consider the following presentation of the first class of all
$n+1$-dimensional filiform Leibniz algebras:
\\ $\emph{U}\cup \emph{F}$, where $\emph{U}= \{L(\alpha):
\alpha_{3}(\alpha_{4}+2\alpha_{3}^{2}) \neq 0\}$, $\emph{F}=
\{L(\alpha): \alpha_{3}(\alpha_{4}+2\alpha_{3}^{2}) = 0\}.$

Our main interest will be the cases of open sets i.e. "generic
algebras" cases.

{\bf Theorem 3.} $i)$ Two algebras $L(\alpha)$ and $L(\alpha')$
from $\emph{U}$ are isomorphic if and only if
$$
\rho_{i}(\frac{2\alpha_{3}}{\alpha_{4}+ 2
\alpha_{3}^{2}},\frac{\alpha_{4}}{2 \alpha_{3}^{2}};\alpha) =
\rho_{i}(\frac{2\alpha_{3}'}{\alpha_{4}'+ 2
\alpha_{3}'^{2}},\frac{\alpha_{4}'}{2 \alpha_{3}'^{2}};\alpha')
$$
whenever \ \ $\emph{i}=\overline{3,n-1}.$

$ii)$ For any $(a_3,a_4,...,a_{n-1})\in C^{n-3}$ there is an
algebra $L(\alpha)$ from $\emph{U}$ such that
$$\rho_{i}(\frac{2\alpha_{3}}{\alpha_{4}+ 2
\alpha_{3}^{2}},\frac{\alpha_{4}}{2 \alpha_{3}^{2}};\alpha)= a_i \
\ \mbox{for all} \ \ \emph{i}=\overline{3,n-1}$$

{\bf Proof.} $i).$ Part "if". Let two algebras $L(\alpha)$ and
$L(\alpha')$ be isomorphic that is there exist $A,B \in
\textbf{C}$ such that $A(A+B)\neq 0$ and
$\alpha'=\rho(\frac{1}{A},\frac{B}{A};\alpha)$. Consider algebra
$L(\alpha^{0})$, where
$\alpha^{0}=\rho(\frac{1}{A_{0}},\frac{B_{0}}{A_{0}};\alpha)$ and
$A_{0}=\frac{\alpha_{4}+2\alpha_{3}^{2}}{2\alpha_{3}},$
$B_{0}=\frac{\alpha_{4}(\alpha_{4}+2\alpha_{3}^{2})}{4\alpha_{3}^{3}}.$
So $\alpha=\rho(A,\frac{-B}{A+B};\alpha')$ and
$$\alpha^{0}=\rho(\frac{2\alpha_{3}}{\alpha_{4}+2\alpha_{3}^{2}},
\frac{\alpha_{4}}{2\alpha_{3}^{2}};\alpha)=\rho(\frac{1}{A_{0}},\frac{B_{0}}{A_{0}};
\rho(A,\frac{-B}{A+B};\alpha'))=\rho(\frac{A}{A_{0}},\frac{B_{0}A-A_{0}B}{A(A+B)};\alpha').$$
It is easy to check that
$\frac{A}{A_{0}}=\frac{2\alpha_{3}'}{\alpha_{4}'+2\alpha_{3}'^{2}}$
and $\frac{B_{0}A-A_{0}B}{A(A+B)}=\frac{\alpha_{4}'}{2
\alpha_{3}'^{2}}.$ Therefore $$
\rho(\frac{2\alpha_{3}}{\alpha_{4}+ 2
\alpha_{3}^{2}},\frac{\alpha_{4}}{2 \alpha_{3}^{2}};\alpha) =
\rho(\frac{2\alpha_{3}'}{\alpha_{4}'+ 2
\alpha_{3}'^{2}},\frac{\alpha_{4}'}{2 \alpha_{3}'^{2}};\alpha')
$$ and, in particular, $$
\rho_{i}(\frac{2\alpha_{3}}{\alpha_{4}+ 2
\alpha_{3}^{2}},\frac{\alpha_{4}}{2 \alpha_{3}^{2}};\alpha) =
\rho_{i}(\frac{2\alpha_{3}'}{\alpha_{4}'+ 2
\alpha_{3}'^{2}},\frac{\alpha_{4}'}{2 \alpha_{3}'^{2}};\alpha')
,$$ for all $i=\overline{3,n-1}.$ This procedure can be shown
schematically by the following picture:

 $\qquad \qquad  \qquad
\qquad \qquad  \alpha \qquad \buildrel(A_{0},B_{0})\over
\longrightarrow  \ \quad  \alpha^{0}$

 $ \qquad \qquad  \qquad \qquad  (A,B) \searrow
\qquad \qquad \nearrow (A_0A^{-1},\frac{B_0A-A_0B}{A_0(A+B)})$
\qquad \qquad \qquad \qquad

$ \qquad \qquad \qquad \qquad \qquad
 \qquad \qquad \alpha'$

Part "Only if". Let the equalities \\ $$
\rho_{i}(\frac{2\alpha_{3}}{\alpha_{4}+ 2
\alpha_{3}^{2}},\frac{\alpha_{4}}{2 \alpha_{3}^{2}};\alpha) =
\rho_{i}(\frac{2\alpha_{3}'}{\alpha_{4}'+ 2
\alpha_{3}'^{2}},\frac{\alpha_{4}'}{2 \alpha_{3}'^{2}};\alpha'),\
\ \emph{i}=\overline{3,n-1}
$$
hold. Then it is easy to notice that
$$
\rho_{i}(\frac{2\alpha_{3}}{\alpha_{4}+ 2
\alpha_{3}^{2}},\frac{\alpha_{4}}{2 \alpha_{3}^{2}};\alpha) =
\rho_{i}(\frac{2\alpha_{3}'}{\alpha_{4}'+ 2
\alpha_{3}'^{2}},\frac{\alpha_{4}'}{2 \alpha_{3}'^{2}};\alpha')\
 \mbox{for}\ \emph{i}=\overline{1,2}
$$ as well and therefore
 $\rho(\frac{2\alpha_{3}}{\alpha_{4}+ 2
\alpha_{3}^{2}},\frac{\alpha_{4}}{2 \alpha_{3}^{2}};\alpha)=
\rho(\frac{2\alpha_{3}'}{\alpha_{4}'+ 2
\alpha_{3}'^{2}},\frac{\alpha_{4}'}{2 \alpha_{3}'^{2}};\alpha')$
that means the algebras $L(\alpha)$ and $L(\alpha')$ are
isomorphic to the same algebra and therefore they are isomorphic
to each other.

 Proof of the part $ii)$. The system of equations
$$\rho_i(\frac{2\alpha_{3}}{\alpha_{4}+ 2
\alpha_{3}^{2}},\frac{\alpha_{4}}{2 \alpha_{3}^{2}};\alpha)= a_i \
\ 3\leq i \leq n-1,$$ where $(a_3,a_4,...,a_{n-1})$ is given and
$\alpha= (\alpha_3,\alpha_4,...,\alpha_{n-1},\theta)$ is unknown,
has a solution as far as for any $3\leq i \leq n-1$ in
$\rho_i(\frac{2\alpha_{3}}{\alpha_{4}+ 2
\alpha_{3}^{2}},\frac{\alpha_{4}}{2 \alpha_{3}^{2}};\alpha)$
($\rho_{n-1}(\frac{2\alpha_{3}}{\alpha_{4}+ 2
\alpha_{3}^{2}},\frac{\alpha_{4}}{2 \alpha_{3}^{2}};\alpha)$) only
variables  $\alpha_3,\alpha_4,...,\alpha_{i}$ (respectively,
$\alpha_3,\alpha_4,...,\alpha_{n-1},\theta$ ) occur and each of
these equations is a linear equation with respect to the last
variable occurred in it. This is the end of the proof.

Here are the corresponding list of invariants for some low
dimensional cases. Case of $n=4$ i.e. \textbf{dim L=5:}
$$ \begin{array}{lll} \left\{
\begin{array}{lll}
\rho_3(\frac{2\alpha_{3}}{\alpha_{4}+ 2
\alpha_{3}^{2}},\frac{\alpha_{4}}{2
\alpha_{3}^{2}};\alpha)=(\frac{2\alpha_{3}}{\alpha_{4}+ 2
\alpha_{3}^{2}})^{2}(\theta-\alpha_{4}).
\end{array}
\right.
\end{array}
$$

Case of $n=5$ i.e. \textbf{dim L=6:}
$$ \begin{array}{lll} \left\{
\begin{array}{lll}
\rho_3(\frac{2\alpha_{3}}{\alpha_{4}+ 2
\alpha_{3}^{2}},\frac{\alpha_{4}}{2
\alpha_{3}^{2}};\alpha)=(\frac{2\alpha_{3}}{\alpha_{4}+ 2
\alpha_{3}^{2}})^2
\frac{\alpha_5+5\alpha_{3}\alpha_4+5\alpha_3^3}{\alpha_{3}}-5. \\
\\
\rho_4(\frac{2\alpha_{3}}{\alpha_{4}+ 2
\alpha_{3}^{2}},\frac{\alpha_{4}}{2
\alpha_{3}^{2}};\alpha)=(\frac{2\alpha_{3}}{\alpha_{4}+ 2
\alpha_{3}^{2}})^3(\theta-\alpha_5)+(\frac{2\alpha_{3}}{\alpha_{4}+
2\alpha_{3}^{2}})^2\frac{\alpha_5+5\alpha_{3}\alpha_4+5\alpha_3^3}{\alpha_{3}}-5.
\end{array}
\right.
\end{array}
$$

Case of $n=6$ i.e. \textbf{dim L=7:}
$$ \begin{array}{lll} \left\{
\begin{array}{lll}
\rho_3(\frac{2\alpha_{3}}{\alpha_{4}+ 2
\alpha_{3}^{2}},\frac{\alpha_{4}}{2
\alpha_{3}^{2}};\alpha)=(\frac{2\alpha_{3}}{\alpha_{4}+ 2
\alpha_{3}^{2}})^2
\frac{\alpha_5+5\alpha_{3}\alpha_4+5\alpha_3^3}{\alpha_{3}}-5. \\
\\
\rho_4(\frac{2\alpha_{3}}{\alpha_{4}+ 2
\alpha_{3}^{2}},\frac{\alpha_{4}}{2
\alpha_{3}^{2}};\alpha)=(\frac{2\alpha_{3}}{\alpha_{4}+2\alpha_{3}^{2}})^3
\frac{\alpha_{6}+6\alpha_{3}\alpha_{5}+21\alpha_{3}^2\alpha_{4}+
3\alpha_{4}^{2}+14\alpha_{3}^4}{\alpha_{3}}\\
\qquad \qquad \qquad \qquad \qquad \qquad \qquad \qquad
-(\frac{2\alpha_{3}}{\alpha_{4}+2\alpha_{3}^{2}})^2
\frac{6\alpha_{3}\alpha_{5}+42\alpha_{3}^2\alpha_{4}+3\alpha_{4}^2+
42\alpha_{3}^4}{\alpha_{3}^{2}}+28.\\
\\
\rho_5(\frac{2\alpha_{3}}{\alpha_{4}+ 2
\alpha_{3}^{2}},\frac{\alpha_{4}}{2
\alpha_{3}^{2}};\alpha)=(\frac{2\alpha_{3}}{\alpha_4+2\alpha_{3}^{2}})^4(\theta-\alpha_6)
+(\frac{2\alpha_{3}}{\alpha_4+2\alpha_{3}^{2}})^3(\frac{\alpha_{6}+6\alpha_3\alpha_5+
21\alpha_3^2\alpha_4+3\alpha_4^2+14\alpha_{3}^4}{\alpha_3})\\
\qquad \qquad \qquad \qquad \qquad \qquad \qquad \qquad
-(\frac{2\alpha_{3}}{\alpha_{4}+ 2 \alpha_{3}^{2}})^2
\frac{6\alpha_{3}\alpha_{5}+42\alpha_{3}^2\alpha_{4}+3\alpha_{4}^2+42\alpha_{3}^4}
{\alpha_{3}^2}+28.
\end{array}
\right.
\end{array}
$$
Case of $n=7$ i.e. \textbf{dim L=8:}
$$ \begin{array}{lll} \left\{
\begin{array}{lll}
\rho_3(\frac{2\alpha_{3}}{\alpha_{4}+ 2
\alpha_{3}^{2}},\frac{\alpha_{4}}{2
\alpha_{3}^{2}};\alpha)=(\frac{2\alpha_{3}}{\alpha_{4}+2\alpha_{3}^{2}})^2
\frac{\alpha_{5}+5\alpha_{3}\alpha_{4}+5\alpha_{3}^{3}}{\alpha_{3}}-5. \\
\\
\rho_4(\frac{2\alpha_{3}}{\alpha_{4}+ 2
\alpha_{3}^{2}},\frac{\alpha_{4}}{2
\alpha_{3}^{2}};\alpha)=(\frac{2\alpha_{3}}{\alpha_{4}+2\alpha_{3}^{2}})^3
\frac{\alpha_{6}+6\alpha_{3}\alpha_{5}+21\alpha_{3}^2\alpha_{4}+
3\alpha_{4}^{2}+14\alpha_{3}^4}{\alpha_{3}}\\
\qquad \qquad \qquad \qquad \qquad \qquad
-(\frac{2\alpha_{3}}{\alpha_{4}+2\alpha_{3}^{2}})^2
\frac{6\alpha_{3}\alpha_{5}+42\alpha_{3}^2\alpha_{4}+3\alpha_{4}^2+
42\alpha_{3}^4}{\alpha_{3}^{2}}+28.\\
\\
\rho_5(\frac{2\alpha_{3}}{\alpha_{4}+ 2
\alpha_{3}^{2}},\frac{\alpha_{4}}{2
\alpha_{3}^{2}};\alpha)=(\frac{2\alpha_{3}}{\alpha_{4}+2\alpha_{3}^{2}})^4
\frac{\alpha_{7}+7\alpha_{3}\alpha_{6}+28\alpha_{3}\alpha_{4}^2+
28\alpha_{3}^{2}\alpha_{5}+7\alpha_{4}\alpha_{5}+
84\alpha_{3}^3\alpha_{4}+42\alpha_{3}^{5}}{\alpha_{3}}\\
\qquad \qquad \qquad \qquad
-(\frac{2\alpha_{3}}{\alpha_{4}+2\alpha_{3}^{2}})^3
\frac{7\alpha_{3}\alpha_{6}+56\alpha_{3}\alpha_{4}^2+56\alpha_{3}^{2}\alpha_{5}+7\alpha_{4}\alpha_{5}+
252\alpha_{3}^3\alpha_{4}+168\alpha_{3}^{5}}{\alpha_{3}^2}\\
\qquad \qquad \qquad \qquad \qquad \qquad \qquad \quad
+(\frac{2\alpha_{3}}{\alpha_{4}+2\alpha_{3}^{2}})^2\frac{28(\alpha_{4}^2+\alpha_{3}\alpha_{5}
+9\alpha_{3}^2\alpha_{4}+9\alpha_{3}^{4})}{\alpha_{3}^2}-126.\\
\\
\rho_6(\frac{2\alpha_{3}}{\alpha_{4}+ 2
\alpha_{3}^{2}},\frac{\alpha_{4}}{2
\alpha_{3}^{2}};\alpha)=(\frac{2\alpha_{3}}{\alpha_{4}+2\alpha_{3}^{2}})^5(\theta-\alpha_{7}) \\
\qquad \qquad \qquad
+(\frac{2\alpha_{3}}{\alpha_{4}+2\alpha_{3}^{2}})^4
\frac{\alpha_{7}+7\alpha_{3}\alpha_{6}+28\alpha_{3}\alpha_{4}^2+
28\alpha_{3}^{2}\alpha_{5}+7\alpha_{4}\alpha_{5}+
84\alpha_{3}^3\alpha_{4}+42\alpha_{3}^{5}}{\alpha_{3}}\\
\qquad \qquad \qquad \qquad
-(\frac{2\alpha_{3}}{\alpha_{4}+2\alpha_{3}^{2}})^3
\frac{7\alpha_{3}\alpha_{6}+56\alpha_{3}\alpha_{4}^2+56\alpha_{3}^{2}\alpha_{5}+7\alpha_{4}\alpha_{5}+
252\alpha_{3}^3\alpha_{4}+168\alpha_{3}^{5}}{\alpha_{3}^2}\\
\qquad \qquad \qquad \qquad \qquad \qquad \qquad \quad
+(\frac{2\alpha_{3}}{\alpha_{4}+2\alpha_{3}^{2}})^2\frac{28(\alpha_{4}^2+\alpha_{3}\alpha_{5}
+9\alpha_{3}^2\alpha_{4}+9\alpha_{3}^{4})}{\alpha_{3}^2}-126.
\end{array}
\right.
\end{array}
$$

The following two theorems deal with the isomorphism problem for
elements from the closed set F. The proof of these theorems are
similar to the proof of Theorem 3. The set $\emph{F}$ in it's turn
can be represented as a union of two open and one closed subsets:
$\emph{F}=\emph{U}'_{1}\bigcup\emph{U}'_{2}\bigcup\emph{F}',$
where $\emph{U}'_{1}=\{L(\alpha)\in\emph{F}: \alpha_{3}\neq0 \ \
\mbox{and} \ \ \alpha_{4}+2\alpha_{3}^{2}= 0\}$,
$\emph{U}'_{2}=\{L(\alpha)\in\emph{F}: \alpha_{3}=0 \ \ \mbox{and}
\ \ \alpha_{4}+2\alpha_{3}^{2}\neq 0\}$, and
$\emph{F}'=\{L(\alpha)\in\emph{F}: \alpha_{3}=0 \ \ \mbox{and} \ \
\alpha_{4}+2\alpha_{3}^{2}=0 \}.$

Then represent $\emph{U}'_{1}$ and $\emph{U}'_{2}$ in the form
$$\emph{U}'_{1}=\emph{U}''_{1}\bigcup\emph{F}''_{1},$$
where $\emph{U}''_{1}=\{\alpha\in\emph{U}'_{1}:
(\alpha_{5}-5\alpha_{3}^{3})(\alpha_{6}+6\alpha_{3}\alpha_{5}-16\alpha_{3}^{4})\neq0\},$
and $\emph{F}''_{1}=\{\alpha\in\emph{U}'_{1}:
(\alpha_{5}-5\alpha_{3}^{3})(\alpha_{6}+6\alpha_{3}\alpha_{5}-16\alpha_{3}^{4})=0\},$
and
$$\emph{U}'_{2}=\emph{U}''_{2}\bigcup\emph{F}''_{2},$$
where
$$\emph{U}''_{2}=\{\alpha\in\emph{U}'_{2}:
\alpha_{5}\neq0\},$$
$$\emph{F}''_{2}=\{\alpha\in\emph{U}'_{2}:
\alpha_{5}=0\},$$

{\bf Theorem 4.} Let $n\geq6.$ Then two algebras $L(\alpha)$ and
$L(\alpha')$ from $\emph{U}''_{1}$ are isomorphic if and only if

$$
\rho_{i}(\frac{5\alpha_{3}^{3}-\alpha_{5}}{\alpha_{6}+6\alpha_{3}\alpha_{5}
-16\alpha_{3}^{4}},\frac{\alpha_{6}+7\alpha_{3}\alpha_{5}-21\alpha_{3}^{4}}
{\alpha_{3}(5\alpha_{3}^{3}-\alpha_{5})};\alpha)=
\rho_{i}(\frac{5\alpha_{3}'^{3}-\alpha_{5}'}{\alpha_{6}'+6\alpha_{3}'\alpha_{5}'
-16\alpha_{3}'^{4}},\frac{\alpha_{6}'+7\alpha_{3}'\alpha_{5}'-21\alpha_{3}'^{4}}
{\alpha_{3}'(5\alpha_{3}'^{3}-\alpha_{5}')};\alpha')$$ for $\ \
\overline{i=4,n-1}.$ Moreover,
$$\rho_{4}(\frac{5\alpha_{3}^{3}-\alpha_{5}}{\alpha_{6}+6\alpha_{3}\alpha_{5}
-16\alpha_{3}^{4}},\frac{\alpha_{6}+7\alpha_{3}\alpha_{5}-21\alpha_{3}^{4}}
{\alpha_{3}(5\alpha_{3}^{3}-\alpha_{5})};\alpha)\neq-14.
$$
and for any $(a_4,...,a_{n-1})\in C^{n-4}, a_4\neq-14$ there is an
algebra $L(\alpha)$ from $\emph{U}''_{1}$ such that
$$\rho_{i}(\frac{5\alpha_{3}^{3}-\alpha_{5}}{\alpha_{6}+6\alpha_{3}\alpha_{5}
-16\alpha_{3}^{4}},\frac{\alpha_{6}+7\alpha_{3}\alpha_{5}-21\alpha_{3}^{4}}
{\alpha_{3}(5\alpha_{3}^{3}-\alpha_{5})};\alpha)=a_i \ \ \mbox{for
all} \ \ \emph{i}=\overline{4,n-1}$$

Here are the list of invariants for $n=6$ and $n=7$ cases:
\\ $\textbf{n=6:}$
$$ \begin{array}{lll} \left\{
\begin{array}{lll}
\rho_{4}(\frac{5\alpha_{3}^{3}-\alpha_{5}}{\alpha_{6}+6\alpha_{3}\alpha_{5}
-16\alpha_{3}^{4}},\frac{\alpha_{6}+7\alpha_{3}\alpha_{5}-21\alpha_{3}^{4}}
{\alpha_{3}(5\alpha_{3}^{3}-\alpha_{5})};\alpha)=7\frac{(5\alpha_{3}^{3}-\alpha_{5})^3}
{\alpha_3(\alpha_{6}+6\alpha_{3}\alpha_{5}
-16\alpha_{3}^{4})^2}-14.\\
\\
\rho_{5}(\frac{5\alpha_{3}^{3}-\alpha_{5}}{\alpha_{6}+6\alpha_{3}\alpha_{5}
-16\alpha_{3}^{4}},\frac{\alpha_{6}+7\alpha_{3}\alpha_{5}-21\alpha_{3}^{4}}
{\alpha_{3}(5\alpha_{3}^{3}-\alpha_{5})};\alpha)=(\frac{5\alpha_{3}^{3}-\alpha_{5}}
{\alpha_{6}+6\alpha_{3}\alpha_{5}
-16\alpha_{3}^{4}})^4(\theta-\alpha_6)+7\frac{(5\alpha_{3}^{3}-\alpha_{5})^3}
{\alpha_3(\alpha_{6}+6\alpha_{3}\alpha_{5}
-16\alpha_{3}^{4})^2}-14.\\
\end{array}
\right.
\end{array}
$$
$\textbf{n=7:}$

$$
\begin{array}{l} \left\{
\begin{array}{l}
\rho_{4}(\frac{5\alpha_{3}^{3}-\alpha_{5}}{\alpha_{6}+6\alpha_{3}\alpha_{5}
-16\alpha_{3}^{4}},\frac{\alpha_{6}+7\alpha_{3}\alpha_{5}-21\alpha_{3}^{4}}
{\alpha_{3}(5\alpha_{3}^{3}-\alpha_{5})};\alpha)=7\frac{(5\alpha_{3}^{3}-\alpha_{5})^3}
{\alpha_3(\alpha_{6}+6\alpha_{3}\alpha_{5}
-16\alpha_{3}^{4})^2}-14.\\
\\
\rho_{5}(\frac{5\alpha_{3}^{3}-\alpha_{5}}{\alpha_{6}+6\alpha_{3}\alpha_{5}
-16\alpha_{3}^{4}},\frac{\alpha_{6}+7\alpha_{3}\alpha_{5}-21\alpha_{3}^{4}}
{\alpha_{3}(5\alpha_{3}^{3}-\alpha_{5})};\alpha)=(\frac{5\alpha_{3}^{3}-\alpha_{5}}
{\alpha_{6}+6\alpha_{3}\alpha_{5}
-16\alpha_{3}^{4}})^4\frac{\alpha_7+7\alpha_{3}\alpha_6-14\alpha_3^2\alpha_{5}-14\alpha_3^5}
{\alpha_3}
\\\qquad \qquad \qquad \qquad \qquad \qquad \qquad \qquad \qquad \qquad -35\frac{(5\alpha_{3}^{3}-\alpha_{5})^3}{\alpha_3(\alpha_{6}+6\alpha_{3}\alpha_{5}
-16\alpha_{3}^{4})^2}+42.\\
\\
\rho_{6}(\frac{5\alpha_{3}^{3}-\alpha_{5}}{\alpha_{6}+6\alpha_{3}\alpha_{5}
-16\alpha_{3}^{4}},\frac{\alpha_{6}+7\alpha_{3}\alpha_{5}-21\alpha_{3}^{4}}
{\alpha_{3}(5\alpha_{3}^{3}-\alpha_{5})};\alpha)=(\frac{5\alpha_{3}^{3}-\alpha_{5}}
{\alpha_{6}+6\alpha_{3}\alpha_{5}
-16\alpha_{3}^{4}})^5(\theta-\alpha_7)\\
\qquad \qquad +(\frac{5\alpha_{3}^{3}-\alpha_{5}}
{\alpha_{6}+6\alpha_{3}\alpha_{5}
-16\alpha_{3}^{4}})^4\frac{\alpha_7+7\alpha_{3}\alpha_6-14\alpha_3^2\alpha_{5}-14\alpha_3^5}
{\alpha_3}-35\frac{(5\alpha_{3}^{3}-\alpha_{5})^3}{\alpha_3(\alpha_{6}+6\alpha_{3}\alpha_{5}
-16\alpha_{3}^{4})^2}+42.
\end{array}
\right.
\end{array}
$$

{\bf Remark.} An analog of Theorem 4 can be stated for $n=4$ and
$n=5$ cases as well. But in these cases it is specific and
therefore is not considered here.

{\bf Theorem 5.} Let $n\geq5$. Then two algebras $L(\alpha)$ and
$L(\alpha')$ from $\emph{U}''_{2}$ are isomorphic if and only if
$$
\rho_{i}(\frac{\alpha_{4}}{\alpha_{5}},\frac{\alpha_{5}^{2}}
{\alpha_{3}^{3}}-1;\alpha) =
\rho_{i}(\frac{\alpha_{4}'}{\alpha_{5}'},\frac{\alpha_{5}'^{2}}
{\alpha_{3}'^{3}}-1;\alpha')$$ for $\ \ \overline{i=4,n-1},$

Moreover, for any $(a_4,...,a_{n-1})\in C^{n-4}$  there is an
algebra $L(\alpha)$ from $\emph{U}''_{2}$ such that
$$\rho_{i}(\frac{\alpha_{4}}{\alpha_{5}},\frac{\alpha_{5}^{2}}
{\alpha_{3}^{3}}-1;\alpha)=a_{i}$$ for all $i=\overline{4,n-1}.$

Here are the corresponding list of invariants for $n=5,6$ and
$n=7$ cases:

Case of $n=5$ i.e. \textbf{dim L=6:}
$$ \begin{array}{lll} \left\{
\begin{array}{lll}
\rho_4(\frac{\alpha_{4}}{\alpha_{5}},\frac{\alpha_{5}^{2}-\alpha_{4}^{3}}{
\alpha_{4}^{3}};\alpha)=\frac{\alpha_{4}(\alpha_4^2\theta-\alpha_{4}^{3}\alpha_3
-3\alpha_{5}^{3})}{\alpha_{5}^{3}}.
\end{array}
\right.
\end{array}
$$

Case of $n=6$ i.e. \textbf{dim L=7:}
$$ \begin{array}{lll} \left\{
\begin{array}{lll}
\rho_4(\frac{\alpha_{4}}{\alpha_{5}},\frac{\alpha_{5}^{2}-\alpha_{4}^{3}}{
\alpha_{4}^{3}};\alpha)=\frac{\alpha_{4}(\alpha_{6}+3\alpha_4^2)}{\alpha_{5}^{2}}
-3; \\ \\
\rho_5(\frac{\alpha_{4}}{\alpha_{5}},\frac{\alpha_{5}^{2}-\alpha_{4}^{3}}{
\alpha_{4}^{3}};\alpha)=(\frac{\alpha_{4}}{\alpha_{5}})^4(\theta-\alpha_{6})+
\frac{\alpha_{4}(\alpha_{6}+3\alpha_4^2)}{\alpha_{5}^{2}} -3.
\end{array}
\right.
\end{array}
$$

Case of $n=7$ i.e. \textbf{dim L=8:}
$$ \begin{array}{lll} \left\{
\begin{array}{lll}

\rho_4(\frac{\alpha_{4}}{\alpha_{5}},\frac{\alpha_{5}^{2}-\alpha_{4}^{3}}{
\alpha_{4}^{3}};\alpha)=\frac{\alpha_{4}(\alpha_{6}+3\alpha_4^2)}{\alpha_{5}^{2}}
-3; \\ \\
\rho_5(\frac{\alpha_{4}}{\alpha_{5}},\frac{\alpha_{5}^{2}-\alpha_{4}^{3}}{
\alpha_{4}^{3}};\alpha)=\frac{\alpha_{4}^{2}(\alpha_{7}+7\alpha_4\alpha_{5})}{\alpha_{5}^{3}}-7; \\
\\
\rho_6(\frac{\alpha_{4}}{\alpha_{5}},\frac{\alpha_{5}^{2}-\alpha_{4}^{3}}{
\alpha_{4}^{3}};\alpha)=(\frac{\alpha_{4}}{\alpha_5})^5(\theta-\alpha_7)
+\frac{\alpha_{4}^{2}(\alpha_{7}+7\alpha_4\alpha_{5})}{\alpha_{5}^{3}}-7.

\end{array}
\right. \\
\\
\end{array}
$$

For any given low dimensional case the above suggested approach
enable us to get the complete classification of filiform Leibniz
algebras from the first class. It is hoped that we can present it
in the near future.

{\bf Acknowledgement.} We would like to thank J.R.Gomez and
B.A.Omirov for their kind permission to use their unpublished
result.

\end{document}